\catcode`\@=11

\font\tenmsx=msxm10
\font\sevenmsx=msxm7
\font\fivemsx=msxm5
\font\tenmsy=msym10
\font\sevenmsy=msym7
\font\fivemsy=msym5
\newfam\msxfam
\newfam\msyfam
\textfont\msxfam=\tenmsx  \scriptfont\msxfam=\sevenmsx
  \scriptscriptfont\msxfam=\fivemsx
\textfont\msyfam=\tenmsy  \scriptfont\msyfam=\sevenmsy
  \scriptscriptfont\msyfam=\fivemsy

\def\hexnumber@#1{\ifnum#1<10 \number#1\else
 \ifnum#1=10 A\else\ifnum#1=11 B\else\ifnum#1=12 C\else
 \ifnum#1=13 D\else\ifnum#1=14 E\else\ifnum#1=15 F\fi\fi\fi\fi\fi\fi\fi}

\def\msx@{\hexnumber@\msxfam}
\def\msy@{\hexnumber@\msyfam}
\mathchardef\boxdot="2\msx@00
\mathchardef\boxplus="2\msx@01
\mathchardef\boxtimes="2\msx@02
\mathchardef\square="0\msx@03
\mathchardef\blacksquare="0\msx@04
\mathchardef\centerdot="2\msx@05
\mathchardef\lozenge="0\msx@06
\mathchardef\blacklozenge="0\msx@07
\mathchardef\circlearrowright="3\msx@08
\mathchardef\circlearrowleft="3\msx@09
\mathchardef\rightleftharpoons="3\msx@0A
\mathchardef\leftrightharpoons="3\msx@0B
\mathchardef\boxminus="2\msx@0C
\mathchardef\Vdash="3\msx@0D
\mathchardef\Vvdash="3\msx@0E
\mathchardef\vDash="3\msx@0F
\mathchardef\twoheadrightarrow="3\msx@10
\mathchardef\twoheadleftarrow="3\msx@11
\mathchardef\leftleftarrows="3\msx@12
\mathchardef\rightrightarrows="3\msx@13
\mathchardef\upuparrows="3\msx@14
\mathchardef\downdownarrows="3\msx@15
\mathchardef\upharpoonright="3\msx@16

\mathchardef\downharpoonright="3\msx@17
\mathchardef\upharpoonleft="3\msx@18
\mathchardef\downharpoonleft="3\msx@19
\mathchardef\rightarrowtail="3\msx@1A
\mathchardef\leftarrowtail="3\msx@1B
\mathchardef\leftrightarrows="3\msx@1C
\mathchardef\rightleftarrows="3\msx@1D
\mathchardef\Lsh="3\msx@1E
\mathchardef\Rsh="3\msx@1F
\mathchardef\rightsquigarrow="3\msx@20
\mathchardef\leftrightsquigarrow="3\msx@21
\mathchardef\looparrowleft="3\msx@22
\mathchardef\looparrowright="3\msx@23
\mathchardef\circeq="3\msx@24
\mathchardef\succsim="3\msx@25
\mathchardef\gtrsim="3\msx@26
\mathchardef\gtrapprox="3\msx@27
\mathchardef\multimap="3\msx@28
\mathchardef\therefore="3\msx@29
\mathchardef\because="3\msx@2A
\mathchardef\doteqdot="3\msx@2B

\mathchardef\triangleq="3\msx@2C
\mathchardef\precsim="3\msx@2D
\mathchardef\lesssim="3\msx@2E
\mathchardef\lessapprox="3\msx@2F
\mathchardef\eqslantless="3\msx@30
\mathchardef\eqslantgtr="3\msx@31
\mathchardef\curlyeqprec="3\msx@32
\mathchardef\curlyeqsucc="3\msx@33
\mathchardef\preccurlyeq="3\msx@34
\mathchardef\leqq="3\msx@35
\mathchardef\leqslant="3\msx@36
\mathchardef\lessgtr="3\msx@37
\mathchardef\backprime="0\msx@38
\mathchardef\risingdotseq="3\msx@3A
\mathchardef\fallingdotseq="3\msx@3B
\mathchardef\succcurlyeq="3\msx@3C
\mathchardef\geqq="3\msx@3D
\mathchardef\geqslant="3\msx@3E
\mathchardef\gtrless="3\msx@3F
\mathchardef\sqsubset="3\msx@40
\mathchardef\sqsupset="3\msx@41
\mathchardef\trianglerighteq="3\msx@44
\mathchardef\trianglelefteq="3\msx@45
\mathchardef\bigstar="0\msx@46
\mathchardef\between="3\msx@47
\mathchardef\blacktriangledown="0\msx@48
\mathchardef\blacktriangleright="3\msx@49
\mathchardef\blacktriangleleft="3\msx@4A
\mathchardef\blacktriangle="0\msx@4E
\mathchardef\triangledown="0\msx@4F
\mathchardef\eqcirc="3\msx@50
\mathchardef\lesseqgtr="3\msx@51
\mathchardef\gtreqless="3\msx@52
\mathchardef\lesseqqgtr="3\msx@53
\mathchardef\gtreqqless="3\msx@54
\mathchardef\Rrightarrow="3\msx@56
\mathchardef\Lleftarrow="3\msx@57
\mathchardef\veebar="2\msx@59
\mathchardef\barwedge="2\msx@5A
\mathchardef\doublebarwedge="2\msx@5B
\mathchardef\angle="0\msx@5C
\mathchardef\measuredangle="0\msx@5D
\mathchardef\sphericalangle="0\msx@5E
\mathchardef\varpropto="3\msx@5F
\mathchardef\smallsmile="3\msx@60
\mathchardef\smallfrown="3\msx@61
\mathchardef\Subset="3\msx@62
\mathchardef\Supset="3\msx@63
\mathchardef\Cup="2\msx@64

\mathchardef\Cap="2\msx@65

\mathchardef\curlywedge="2\msx@66
\mathchardef\curlyvee="2\msx@67
\mathchardef\leftthreetimes="2\msx@68
\mathchardef\rightthreetimes="2\msx@69
\mathchardef\subseteqq="3\msx@6A
\mathchardef\supseteqq="3\msx@6B
\mathchardef\bumpeq="3\msx@6C
\mathchardef\Bumpeq="3\msx@6D
\mathchardef\lll="3\msx@6E

\mathchardef\ggg="3\msx@6F

\mathchardef\circledS="0\msx@73
\mathchardef\pitchfork="3\msx@74
\mathchardef\dotplus="2\msx@75
\mathchardef\backsim="3\msx@76
\mathchardef\backsimeq="3\msx@77
\mathchardef\complement="0\msx@7B
\mathchardef\intercal="2\msx@7C
\mathchardef\circledcirc="2\msx@7D
\mathchardef\circledast="2\msx@7E
\mathchardef\circleddash="2\msx@7F
\def\ulcorner{\delimiter"4\msx@70\msx@70 }
\def\urcorner{\delimiter"5\msx@71\msx@71 }
\def\llcorner{\delimiter"4\msx@78\msx@78 }
\def\lrcorner{\delimiter"5\msx@79\msx@79 }
\def\yen{\mathhexbox\msx@55 }
\def\checkmark{\mathhexbox\msx@58 }
\def\circledR{\mathhexbox\msx@72 }
\def\maltese{\mathhexbox\msx@7A }
\mathchardef\lvertneqq="3\msy@00
\mathchardef\gvertneqq="3\msy@01
\mathchardef\nleq="3\msy@02
\mathchardef\ngeq="3\msy@03
\mathchardef\nless="3\msy@04
\mathchardef\ngtr="3\msy@05
\mathchardef\nprec="3\msy@06
\mathchardef\nsucc="3\msy@07
\mathchardef\lneqq="3\msy@08
\mathchardef\gneqq="3\msy@09
\mathchardef\nleqslant="3\msy@0A
\mathchardef\ngeqslant="3\msy@0B
\mathchardef\lneq="3\msy@0C
\mathchardef\gneq="3\msy@0D
\mathchardef\npreceq="3\msy@0E
\mathchardef\nsucceq="3\msy@0F
\mathchardef\precnsim="3\msy@10
\mathchardef\succnsim="3\msy@11
\mathchardef\lnsim="3\msy@12
\mathchardef\gnsim="3\msy@13
\mathchardef\nleqq="3\msy@14
\mathchardef\ngeqq="3\msy@15
\mathchardef\precneqq="3\msy@16
\mathchardef\succneqq="3\msy@17
\mathchardef\precnapprox="3\msy@18
\mathchardef\succnapprox="3\msy@19
\mathchardef\lnapprox="3\msy@1A
\mathchardef\gnapprox="3\msy@1B
\mathchardef\nsim="3\msy@1C
\mathchardef\napprox="3\msy@1D
\mathchardef\nsubseteqq="3\msy@22
\mathchardef\nsupseteqq="3\msy@23
\mathchardef\subsetneqq="3\msy@24
\mathchardef\supsetneqq="3\msy@25
\mathchardef\subsetneq="3\msy@28
\mathchardef\supsetneq="3\msy@29
\mathchardef\nsubseteq="3\msy@2A
\mathchardef\nsupseteq="3\msy@2B
\mathchardef\nparallel="3\msy@2C
\mathchardef\nmid="3\msy@2D
\mathchardef\nshortmid="3\msy@2E
\mathchardef\nshortparallel="3\msy@2F
\mathchardef\nvdash="3\msy@30
\mathchardef\nVdash="3\msy@31
\mathchardef\nvDash="3\msy@32
\mathchardef\nVDash="3\msy@33
\mathchardef\ntrianglerighteq="3\msy@34
\mathchardef\ntrianglelefteq="3\msy@35
\mathchardef\ntriangleleft="3\msy@36
\mathchardef\ntriangleright="3\msy@37
\mathchardef\nleftarrow="3\msy@38
\mathchardef\nrightarrow="3\msy@39
\mathchardef\nLeftarrow="3\msy@3A
\mathchardef\nRightarrow="3\msy@3B
\mathchardef\nLeftrightarrow="3\msy@3C
\mathchardef\nleftrightarrow="3\msy@3D
\mathchardef\divideontimes="2\msy@3E
\mathchardef\varnothing="0\msy@3F
\mathchardef\nexists="0\msy@40
\mathchardef\mho="0\msy@66
\mathchardef\thorn="0\msy@67
\mathchardef\beth="0\msy@69
\mathchardef\gimel="0\msy@6A
\mathchardef\daleth="0\msy@6B
\mathchardef\lessdot="3\msy@6C
\mathchardef\gtrdot="3\msy@6D
\mathchardef\ltimes="2\msy@6E
\mathchardef\rtimes="2\msy@6F
\mathchardef\shortmid="3\msy@70
\mathchardef\shortparallel="3\msy@71
\mathchardef\smallsetminus="2\msy@72
\mathchardef\thicksim="3\msy@73
\mathchardef\thickapprox="3\msy@74
\mathchardef\approxeq="3\msy@75
\mathchardef\succapprox="3\msy@76
\mathchardef\precapprox="3\msy@77
\mathchardef\curvearrowleft="3\msy@78
\mathchardef\curvearrowright="3\msy@79
\mathchardef\digamma="0\msy@7A
\mathchardef\varkappa="0\msy@7B
\mathchardef\hslash="0\msy@7D
\mathchardef\hbar="0\msy@7E
\mathchardef\backepsilon="3\msy@7F
\def\Bbb{\ifmmode\let\next\Bbb@\else
 \def\next{\errmessage{Use \string\Bbb\space only in math mode}}\fi\next}
\def\Bbb@#1{{\Bbb@@{#1}}}
\def\Bbb@@#1{\fam\msyfam#1}

\catcode`\@=12

\def\C{\Bbb C} 
\def\R{\Bbb R} 
\def\P{\Bbb P} 

\documentstyle{article}

\begin{document}
\sloppy
\newtheorem{Th}{Theorem}[section]
\newtheorem{Satz}{Satz}[section]
\newtheorem{Prop}{Proposition}[section]
\newtheorem{Lemma}{Lemma}[section]
\newtheorem{Rem}{Remark}[section]
\newtheorem{Def}{Definition}[section]
\newtheorem{Cor}{Corollary}[section]
\author{D.~Zaitsev}
\title{On the linearization of the automorphism groups of algebraic domains.}
\maketitle
\section{Linearization Theorem and applications}

	Let $D$ be a domain in $\C^n$ and
$G$ a topological group which acts effectively
on $D$ by holomorphic automorphisms. In this paper we are interested in
projective linearizations of the action of $G$, i.e. a linear representation
of $G$
in some $\C^{N+1}$ and an equivariant imbedding of $D$ into $\P^N$ with
respect to this representation. Since $G$ acts effectively, the representation
in $\C^{N+1}$ must be faithful.
In our previous paper~\cite{Z}, however, we considered
an example of a bounded domain $D\subset\C^2$ with an effective action of a
finite covering $G$ of the group $SL_2(\R)$. In this case the group $G$ doesn't
admit a faithful representation. The example shows that a linearization
in the above sense doesn't exist in general.

	In the present paper we give a criterion for the existence of the
projective linearization for birational automorphisms.
The domains we discuss here are open connected sets defined by finitely
many real polynomial inequalities or connected finite unions of such sets.
These domains are called {\it algebraic}.
For instance, in the above example the domain $D$ is algebraic.

\begin{Def}\label{Nash}
\begin{enumerate}
\item A {\bf Nash map} is a real analytic map
$$f=(f_1,\ldots,f_m)\colon U\to \R^m$$ (where $U\subset\R^n$ is open)
such that for each of the components $f_k$ there is a
non-trivial polynomial $P_k$ with $$P_k(x_1,\ldots,x_n,f_k(x_1,\ldots,x_n))=0$$
for all $(x_1,\ldots,x_n)\in U$.
\item A {\bf Nash manifold} $M$ is a real analytic manifold with finitely many
coordinate charts $\phi_i\colon U_i\to V_i$ such that $V_i\subset\R^n$ is
Nash for all $i$ and the transition functions are Nash
(a Nash atlas).
\item A {\bf Nash group} is a Nash manifold with a group operation
$(x,y)\mapsto xy^{-1}$ which is Nash with respect to all Nash coordinate
charts.
\end{enumerate}
\end{Def}

	In the above example the group $SL_2(\R)$ and its finite covering $G$
are Nash groups (The universal covering of $SL_2(\R)$ is a so-called
locally Nash group).
Moreover, the action $G\times D\to D$ is also Nash.
Since the linearization doesn't exist here,
we need a stronger condition on the action of $G$.

	A topological group $G$ is said to be
a {\it group of birational automorphisms}
of a domain $D\subset\C^n$ if we are given an effective (continuous)
action $G\times D\to D$ such that every element
$g\in G$ defines an automorphism of $D$ which extends to a birational
automorphism of $\C^n$. By the {\it degree} of a Nash map $f$ we mean the
minimal natural number $d$ such that all polynomials $P_k$
in Definition~\ref{Nash} can be chosen such that their degrees don't
exceed $d$.
Finally, under a {\it biregular} map between two algebraic
varieties we understand an isomorphism in sense of algebraic geometry.
Our main result is the following linearization criterion.
It will be proved in section~\ref{proof}.

\begin{Th}\label{main}
  Let  $D\subset\C^n$ be a algebraic domain and $G$ a group
of birational automorphisms of $D$. The following properties are equivalent:
\begin{enumerate}
\item $G$ is a subgroup of a Lie group $\hat G$ of birational automorphisms
of $D$ which extends the action of $G$ and has finitely many connected
components;
\item $G$ is a subgroup of a Nash group $\hat G$ of birational
automorphisms of $D$ which extends the action of $G$ to a Nash action
$\hat G\times D\to D$;
\item $G$ is a subgroup of a Nash group $\hat G$ such that the action
$G\times D\to D$ extends to a Nash action $\hat G\times D\to D$;
\item the degree of the automorphism $\phi_g\colon D\to D$
defined by $g\in G$ is bounded;
\item there exists a projective linearization,
i.e. a linear representation of $G$ in some $\C^{N+1}$ and
a biregular imbedding $i\colon \P^n \hookrightarrow \P^N$ such that
the restriction $i|_D$ is $G$-equivariant.
\end{enumerate}
\end{Th}

	We finish this section with applications of Theorem~\ref{main}.
In the previous paper~(\cite{Z}) we gave sufficient conditions on $D$ and $G$
such that $G$ is a Nash group and the action $G\times D\to D$ is Nash.
The condition on $D$ is to be bounded and
to have a non-degenerate boundary in the following sense.

\begin{Def}\label{deg}
	A boundary of a domain $D\subset\C^n$ is called {\bf non-degenerate}
if it contains a smooth point where the Levi-form is non-degenerate.
\end{Def}

	The group $G$ is taken to be the group $Aut_a(D)$ of all
holomorphic Nash (algebraic) automorphisms of $D$. We proved in~\cite{Z} that,
if $D$ is a algebraic bounded domain with non-degenerate boundary,
the group $Aut_a(D)$ is closed in $Aut(D)$ and carries
a unique structure of a Nash group such that
the action $Aut_a(D)\times D\to D$ is Nash with respect to this structure.

	Now let $G=Aut_b(D)\subset Aut_a(D)$ be the group of all birational
automorphisms of $D$. Then $G$ satisfies the property~3 in Theorem~\ref{main}
with $\hat G = Aut_a(D)$. By the property~2, $G$ is a subgroup of a Nash
group of birational automorphisms of $D$.
Since $G$ contains all the birational automorphisms of $D$, $G$ is itself
a Nash group with the Nash action on $D$. We obtain the following corollary.

\begin{Cor}\label{aut-b}
	Let $D\subset\subset\C^n$ be a bounded algebraic domain with non-degenerate
boundary.
Then the group $Aut_b(D)$ of all birational automorphisms of
$D$ is Nash with the Nash action on $D$ which
admits a projective linearization, i.e. there exist a representation
of $Aut_b(D)$ in some $\C^{N+1}$ and
a biregular imbedding $i\colon \P^n \hookrightarrow \P^N$ such that
the restriction $i|_D$ is $Aut_b(D)$-equivariant.
\end{Cor}

	Furthermore, S.~Webster (see \cite{W}) established the following
sufficient conditions on $D$ which make its automorphisms birational.
Let $D$ be a algebraic domain. The theory of semialgebraic sets
(see Benedetti-Risler~\cite{BR}) implies that the boundary
$\partial D$ is contained in finitely many irreducible real hypersurfaces.
Several of them, let say $M_1,\ldots,M_k$, have generically
non-degenerate Levi forms. If $\partial D$ is non-degenerate in sense of
Definition~\ref{deg},
such hypersurfaces exist. The complexifications ${\cal M}_i$'s
of $M_i$'s are defined
to be their complex Zariski closures in $\C^n\times\overline{\C^n}$ where
$M_i$'s are totally real imbedded via the diagonal map $z\mapsto (z,\bar z)$.
It follows that ${\cal M}_i$'s are the irreducible complex hypersurfaces.
Furthermore, the so-called Segre varieties $Q_{iw}$'s,
$w\in \C^n$ are defined by
$$ Q_{iw} := \{z\in \C^n \mid (z, \bar w) \in {\cal M}_i \}.$$
The complexifications and Segre varieties are the important biholomorphic
invariants of a domain $D$ and play a decisive role in the reflection
principle.

	Now we are ready to formulate the conditions of S.~Webster.
\begin{Def}
	A algebraic domain is said to satisfy the condition $(W)$
if for all $i$ the Segre varieties
$Q_{iw}$ uniquely determine $z\in\C^n$ and $Q_{iw}$ is an irreducible
hypersurface in $\C^n$ for all $z$ off a proper subvariety $V_i\subset\C^n$.
\end{Def}

	The Theorem of S.~Webster (see \cite{W}, Theorem~3.5) can be formulated
in the following form:

\begin{Th}\label{bir}
	Let $D\subset\C^n$ be a algebraic domain with non-degenerate boundary
which satisfies the condition $(W)$. Further,
let $f\in Aut(D)$ be an automorphism which is holomorphically
extendible to a smooth boundary point
with non-degenerate Levi-form. Then $f$ is birationally extendible to the
whole $\C^n$.
\end{Th}

	Since every Nash automorphism $f\in Aut_a(D)$ extends holomorphically
to generic boundary points, we obtain the following Corollary.

\begin{Cor}
	Let $D\subset\C^n$ be a bounded algebraic domain which satisfies
the condition $(W)$. Then the whole
group $Aut_a(D)$ is projective linearizable,
i.e. there exist a representation
of $Aut_a(D)$ in some $\C^{N+1}$ and
a biregular imbedding $i\colon \P^n \hookrightarrow \P^N$ such that
the restriction $i|_D$ is $Aut_a(D)$-equivariant.
\end{Cor}

	To obtain the extendibility of the whole group $Aut(D)$ of holomorphic
automorphisms, we consider
the {\it algebraic} domains in sense of Diederich-Forn\ae ss (see \cite{DF}).

\begin{Def}
	A domain $D\subset\subset C^n$ is called {\bf algebraic} if there exists
a real polynomial $r(z,\bar z)$ such that $D$ is a connected component of the
set $$\{z\in\C^n \mid r(z,\bar z)<0 \}$$
and $dr(z)\ne 0$ for $z\in\partial D$.
\end{Def}

	The following fundamental result for such domains is due to
K.~Diederich and J.~E.~Forn\ae ss (see \cite{DF}).

\begin{Th}\label{hol}
	Let $D\subset\subset\C^n$ be an algebraic domain.  Then $Aut_a(D)=Aut(D)$.
\end{Th}

	Thus we obtain the linearization of the whole automorphism group $Aut(D)$.

\begin{Th}\label{alg}
	Let $D\subset\subset \C^n$ be an algebraic domain which satisfies
the condition $(W)$. Then the group $Aut(D)$ is projective linearizable,
i.e. there exist a representation
of $Aut(D)$ in some $\C^{N+1}$ and
a biregular imbedding $i\colon \P^n \hookrightarrow \P^N$ such that
the restriction $i|_D$ is $Aut(D)$-equivariant.
\end{Th}

Further corollaries are devoted to the constructions of complexifications.

\section{Complexifications}
	Using the linearization criterion we establish here existences of
complexifications.
To every real Lie group $G$ one can associate its complexification
(see Hochschild~\cite{Ho}) defined as follows.

\begin{Def}\label{cpx}
	Let $G$ be a real Lie group. A complex Lie group $G^{\C}$ together with
a Lie homomorphism $\imath\colon G \to G^{\C}$ is called a
{\bf complexification} of $G$ if for a given Lie homomorphism
$\phi$ from $G$ into a complex Lie group $H$, there exists exactly one
holomorphic Lie homomorphism $\phi^{\C}\colon G^{\C} \to H$ such that
$\phi=\phi^{\C}\circ\imath$.
A real Lie group $G$ is called {\bf holomorphically extendible} if the map
$\imath\colon G\to G^{\C}$ is injective.
\end{Def}

	A complexification always exists and is unique up to biholomorphisms
(see Hochschild~\cite{Ho} and Heinzner~\cite{He} and \cite{He1}).
Further, one defines the complexification of an action
(see Heinzner~\cite{He}).

\begin{Def}\label{G-cpx}
	Let a real Lie group $G$ act on a complex space $X$ by holomorphic
automorphisms. A complex space $X^{\C}$ together with a holomorphic action
of $G^{\C}$ and a $G$-equivariant map $\imath\colon X\to X^{\C}$ is called
a $G$-{\bf complexification} of $X$ if to every holomorphic $G$-equivariant
map $\phi\colon X\to Y$ into another complex space $Y$ with a holomorphic
action of $G^{\C}$ there exists exactly one holomorphic $G^{\C}$-equivariant
map $\phi^{\C}$ such that $\phi=\phi^{\C}\circ\imath$.
\end{Def}

	A $G$-complexification is unique up to biholomorphic $G^{\C}$-equivariant
maps provided it exists. P.~Heinzner proved in \cite{He} the existence of
a $G$-complexification of $X$ with properties that $\imath\colon X\to X^{\C}$
is an open imbedding
and $X^{\C}$ is Stein in case $G$ is compact and $X$ is a Stein space.

	Now the projective linearization in Theorem~\ref{main} implies
the existence of complexifications in our situation.

\begin{Cor}\label{main1}
	Let $D\subset\C^n$ be a algebraic domain and $G$ a Lie group
of birational automorphisms of $D$ which satisfies one of the equivalent
properties in Theorem~\ref{main}. Then the group $G$ is
holomorphically extendible and there exists a smooth $G$-complexification
$D^{\C}$ of $D$ such the map $\imath\colon D\to D^{\C}$ is an open imbedding.
\end{Cor}

	{\bf Proof.} By property~5 in Theorem~\ref{main}, $G$ is a subgroup of the
complex Lie group $GL_N(\C)$. By Definition~\ref{cpx}, $G$ is
holomorphically extendible. Let $i$ be the embedding of
$\C^n\supset D$, given by Theorem~\ref{main}. Since the stability group
$H\subset GL_N(\C)$ of the complex projective variety
$X:=i(\P^n)$ is a complex Lie group and $G\subset H$,
the Definition~\ref{cpx} yields a holomorphic
action of $G^{\C}$ on $X$. We claim that
$D^{\C}:=G^{\C}\cdot D\subset X$ is the required $G$-complexification
of $D$. Indeed, let $\phi\colon D\to Y$ be a $G$-equivariant holomorphic
map into another complex space $Y$ with a holomorphic
action of $G^{\C}$. To define the required in Definition~\ref{G-cpx}
map $\phi^{\C}$ we take a point $z\in D^{\C}$ which is always of the form
$z=Ax$ with $A\in G^{\C}$ and $x\in D$. Then we set
$\phi^{\C}(z):=A\phi(x)$. Why is $\phi^{\C}(z)$ independent of the
representation $z=Ax$? Because the holomorphic map $A\mapsto A\phi(x)$ is
determined by values on the maximal totally real subgroup $G$:
for $A\in Aut_b(D)$ one has $A\phi(x)=\phi(Ax)$. We obtain a well-defined
$G^{\C}$-equivariant map $\phi^{\C}(z)\colon D\to D^{\C}$ with
the property $\phi=\phi^{\C}\circ\imath$ (because for $z\in D$ one can
choose $A=1$). The holomorphicity of $\phi^{\C}$ is obtained by fixing $A$
in the formula $\phi^{\C}(z):=A\phi(x)$. \hfill $\square$

	For the algebraic domains we obtain the following Corollaries.

\begin{Cor}
	Let $D\subset\subset\C^n$ be a algebraic domain with non-degenerate boundary.
Then the group $Aut_b(D)$ is
holomorphically extendible and there exists an $Aut_b(D)$-complexification
of $D$.
\end{Cor}

\begin{Cor}
	Let $D\subset\subset\C^n$ be a algebraic domain with non-degenerate boundary
which satisfies the condition $(W)$.
Then the group $Aut_a(D)$ is
holomorphically extendible and there exists an $Aut_a(D)$-complexification
of $D$.
\end{Cor}

\begin{Cor}
	Let $D\subset\subset \C^n$ be an algebraic domain which satisfies
the condition $(W)$. Then the group $Aut(D)$ is
holomorphically extendible and there exists an $Aut(D)$-complexification
of $D$.
\end{Cor}
\section{Proof of the main Theorem}\label{proof}

	Let $D$ be a algebraic domain and $G$ a group of birational automorphisms
of $D$. We prove the equivalence of the properties in Theorem~\ref{main}
in the direction of the following two chains:
$2\Longrightarrow 3\Longrightarrow 4\Longrightarrow 5\Longrightarrow 2$
and $2\Longrightarrow 1\Longrightarrow 4$.

	\fbox{$2\Longrightarrow 3$.} The proof is trivial. \hfill $\square$
	\bigskip

	\fbox{$3\Longrightarrow 4$.} Let $G$ be a subgroup of
a Nash group $\hat G$ such that the action
$G\times D\to D$ extends to a Nash action $\hat G\times D\to D$.
We prove the statement for arbitrary Nash manifold $\hat G$ and
Nash map $\hat G\times D\to D$
by induction on $\dim G$. It is obvious for $\dim G=0$.
Let $U\subset \hat G$ be a Nash coordinate chart and
$\phi_i(g)\colon D\to \R$ be the $i$th coordinate of $\phi_g\colon D\to D$
for $g\in U$. Since the map $\phi_i\colon U\times D\to \R$ is Nash, it
satisfies a polynomial equation
$P(g,x,\phi_i(g,x))\equiv 0$. This yields polynomial equations of the same
degree for all $g\in U$ outside a proper Nash submanifold. This submanifold
has lower dimension and the statement is true for it by induction.
In summary, we obtain the boundness of the degree for the whole
neighborhood $U$ and, since the Nash atlas is finite, for $G$. \hfill$\square$

	\bigskip
	\fbox{$4\Longrightarrow 5$.} Here is a sketch of the proof.
The idea is to imbed the group
$G$ into a complex algebraic variety so that the action on $D$
is given by a rational mapping. Using this mapping we construct
a collection of homogeneous polynomials on $\C^{n+1}$ which generate
a finite-dimensional linear subspace, invariant with respect to
the action of $G$. These polynomials
yield the required projective linearization.

	The imbedding of $G$ is obtained by associating to every element
$g\in G$ the Chow coordinates of the complex Zariski closure of the graph of
the automorphism defined by $g$ (see Shafarevich, \cite{S} ,page 65).
The main problem here is that the Chow scheme $C$ has infinitely many disjoint
components parameterized by dimensions and degrees of subvarieties.
In order to concern finitely many components of $C$ we have required
the degree of automorphisms $\phi_g\colon D\to D$ to be bounded.

	To every $g\in G$ one associates the $n$-dimensional
(complex) Zariski closure $\tilde\Gamma_g\subset \P^n\times\P^n$ of the
graph $\Gamma_g\subset D\times D$ of the automorphism defined by
$g$.  To regard $\tilde\Gamma_g$ as a
subvariety   of  some  $P^N$  let  us  consider  the  Segre  imbedding:

\begin{equation}\label{segre}
           v([z_0,...,z_n],[w_0,...,w_n])=[{z_iw_j}]_
               {0 \le i \le n, 0 \le j \le n},
\end{equation}

          $$v\colon \P^n\times \P^n \to \P^{n^2+2n}.$$

	We   set   $N=n^2+2n$   and   obtain   a  family  of
subvarieties  $\rho(g):=v(\tilde\Gamma_g)\subset \P^N$ parameterized by $g\in
G$.
The family $V$ of all algebraic subvarieties of $\P^n$ is parameterized
by the Chow scheme $C$ (see Shafarevich, \cite{S} ,page 65).
Different automorphisms $g\in G$ define different subvarieties
$v(\tilde\Gamma_g)\subset \P^N$ and one obtains an imbedding $\rho$of
$G$ in the Chow scheme $C$.
The (complex) dimension of the subvarieties $v(\tilde\Gamma_g)$ is $n$.
The degree of $v(\tilde\Gamma_g)\subset\P^N$ is the
intersection number with $N-n$ generic linear hyperplanes
$\{L_1=0\},\ldots,\{L_{N-n}=0\}$. It is equal to the intersection number of
$\tilde\Gamma_g\subset\P^N\times\P^N$ with divisors $v^*L_1,\ldots,v^*L_{N-n}$.
By the Bezout theorem this intersection number is bounded.

	Thus, $G$ lies in fact in finitely
many components of the Chow scheme $C$. Let $C_0$ denote the union of
these components and $V_0$ the corresponding family of subvarieties of $P^N$.
We obtain an imbedding $\rho$ of $G$ in a complex projective variety $C_0$.

\begin{Lemma}
	The imbedding $\rho\colon G\to C_0$ is continuous.
\end{Lemma}

	{\bf Proof.} Assume the contrary. Then there exists a sequence
$g_n\to g$ in $G$ such that no subsequence $\rho(g_{n(k)})$ converges
to $\rho(g)$. On the other hand, since the degree of $\rho(g_n)$ is bounded,
there exists a subsequence $\rho(g_{n(k)})$ which converges in $C_0$.
This follows from the Theorem of Bishop (see e.g. F.~Campana, \cite{Cam}).
Let $A\in C_0$ be the limes cycle of this subsequence. Our cycles lie
in $v(\P^n\times\P^n)$ and we identify them with the preimages in
$\P^n\times\P^n$. Since the action $G\times D\to D$ is continuous,
the cycle $A$ contains the graph of $\phi_g$ and therefore its Zariski
closure $\rho(g)$.

	We claim that $A=\rho(g)$. This yields a contradiction with the choice
of $\rho(g_{n(k)})$. Indeed, otherwise there exists
a horizontal of vertical $n$-dimensional projective subspace
$H\subset\P^n\times\P^n$ ($H=\{z\}\times\P^n$ or $H=\P^n\times\{w\}$)
such that the intersection number of $A$ and $H$ is more than one.
Then the intersection number of $\rho(g_{n(k)})$ and $H$ is also more than
one which contradicts to the birationality of $\phi_{g_{n(k)}}$. \hfill
$\square$

   Further  let  $\tilde G$  be the complex Zariski closure of $G$ in $C_0$.

\begin{Lemma}
	The action $\phi\colon  G\times D\to D$ extends to a rational map
$\tilde\phi \colon \tilde G\times \P^n \to \P^n$.
\end{Lemma}

	{\bf Proof.} We begin with the construction of the graph
$\Gamma_{\tilde\phi} \subset \tilde G\times \P^n \to \P^n$ of $\tilde\phi$.
For this we regard $\P^n \times \P^n$ as a subset of $P^N$
(via the Segre imbedding $v$ in (\ref{segre})). We then define
$\Phi = \Gamma_{\tilde\phi}$ to be the intersection of the Chow family $V_0$
with $\tilde G\times \P^n \to \P^n$. This is a complex algebraic variety.
Moreover, for $g\in G$ and $x\in D$ the fibre $\Phi_{(g,x)}\subset\P^n$
consists of the single point $g(x)$. Since the set $G\times D$ is
Zariski dense in $\tilde G\times \P^n$, this is true for every generic fibre of
$\Phi$. This means that $\Phi$ is the graph of a rational map
$\tilde\phi \colon \tilde G\times \P^n \to \P^n$. \hfill $\square$

	The projective variety $\tilde G$ is imbedded in a projective space $\P^m$.
The map $\tilde\phi$ can be extended to a rational map from
$\P^m \times \P^n$ into $\P^n$. Such map is given by $n+1$ polynomials
$P_1(x,y),\dots,P_{n+1}(x,y)$, homogeneous separately in
$x\in\C^{m+1}$ and $y\in\C^{n+1}$.
Let $h$ be a fixed homogeneous polynomial on $\C^{n+1}$. Then the function
$$(x,y)\mapsto h(P_1(x,y),\ldots,P_{n+1}(x,y))$$
is a separately homogeneous polynomial on $\C^{m+1} \times \C^{n+1}$.
The algebra $\C_h[x,y]$ of such polynomials is equal to the tensor
product $\C_h[x]\otimes\C_h[y]$. Therefore there exist polynomials
$\varphi_i\in \C_h[x]$, $\psi_i\in \C_h[y]$, $i=1,\ldots,l$ such that
$$h(P_1(x,y),\ldots,P_{n+1}(x,y)) = \sum_{i=1}^l \varphi_i(x) \psi_i(y).$$
For $x=g\in G$ fixed we obtain
$$\alpha_*(f^{-1}) h = \sum_{i=1}^l c_i \psi_i(y),$$
where $\alpha_*$ denotes the associated action of $G$ on homogeneous
polynomials. In other words, the orbit of $h$ via the action of $G$
is contained in the finite-dimensional subspace
$<\psi_1,\ldots,\psi_l> \subset \C_h[y]$. The linear hull of this orbit
is a finite-dimensional $G$-invariant subspace containing $h$.

	We choose mow sufficiently many polynomials $h_j$, $j=1,\ldots,s$ which
separate the points of $\C^{n+1}$ and such that neither $h_j$ nor
the differentials
$dh_j$ nowhere vanish simultaneously. They lie in a finite-dimensional
$G$-invariant subspace $L\subset \C_h[y]$. Let $(p_1,\ldots,p_{N+1})$
be a collection of homogeneous polynomials which yields a basis of $L$.
The required representation of $G$ is the action on $L$ and the polynomial
map $(p_1,\ldots,p_{u+1})\colon \C^{n+1}\to \P^N$
defines the required projective linearization. \hfill $\square$

	\bigskip
	\fbox{$5\Longrightarrow 2$.} Assume we are given a projective linearization
of the action of $G$
on $D$. It follows that the given representation of $G$ is faithful and
we identify $G$ with its image in $GL_{N+1}(\C)$.
We define now the group $\hat G\supset G$ to be the subgroup of all
$g\in GL_{N+1}$ such that $g(i(D))=i(D)$. It follows that $G\subset\hat G$.

	We wish to prove that $\hat G$ is a Nash subgroup of
$GL_{N+1}(\C)$. For the proof we use the technique of
semialgebraic sets and maps which are closely related to the Nash
manifolds and maps.
The semialgebraic subsets of $\R^n$ are the sets of the form
$\{P_1=\cdots=P_k, Q_1<0,\ldots,Q<s\}$ and finite unions of them where
$P_1,\ldots,P_k$ and $Q_1,\ldots,Q_s$ are real polynomials on $\R^n$.
More generally, the semialgebraic subsets of a Nash manifold $M$ are the
subsets which have semialgebraic intersections with every Nash coordinate
chart.
The semialgebraic maps between semialgebraic sets are any maps with
semialgebraic graphs. The Nash submanifolds of $\R^n$ are exactly
semialgebraic real analytic submanifolds and the Nash maps are
semialgebraic real analytic maps.

	Now the graph $\Gamma\subset Gl_{N+1}\times i(D) \to \P^N$ of the restriction
to $i(D)$ of the linear action of $GL_{N+1}$ on $\P^N$ is a semialgebraic
subset. The condition $g(i(D))=i(D)$ on $g$ defines a semialgebraic subset
$\hat G\subset Gl_{N+1}$. We proved this in the previous paper
(see \cite{Z}, Lemma~6.2). Since $\hat G$ is a subgroup, it is Nash.
Thus, $G$ in a subgroup of the Nash group $\hat G$ of birational automorphisms
of $i(D)\cong D$ with required properties.\hfill $\square$

	\bigskip
	\fbox{$2\Longrightarrow 1$.} A Nash open subset of $\R^n$ is semialgebraic
and has therefore finitely many connected components
(see Benedetti-Risler, \cite{BR}, Theorem~2.2.1). The Nash group $\hat G$
admits a finite Nash atlas and has also finitely many components.\hfill
$\square$

	\bigskip
	\fbox{$1\Longrightarrow 4$.} Assume $G$ is a subgroup of a Lie group $\hat G$
of birational automorphisms of $D$ with finitely many connected components.
Consider the complex coordinates $\phi_i\colon \hat G\times D\to \C$ of the
action of $\hat G$. For fixed $g\in \hat G$ the map
$\phi_i(g)\colon D\to \C$ extends to a rational map
$\tilde\phi_i(g)\colon \C^n\to \C$. These extensions define a map
$\tilde\phi_i\colon \hat G\times\C^n\to \C$.
A priori we don't
know  whether  this  new  map  is  real  analytic  or  even continuous.
To prove this we use the following result of Kazaryan (\cite{Ka}):

\begin{Prop}\label{Kaz}
	 Let  $D'$  be  a domain in $\C^n$  and let $E\subset D'$ be a
   nonpluripolar\footnote{a subset $E\subset D'$ is called {\it nonpluripolar}
if
there are no plurisubharmonic functions $f\colon D'\to \R\cup\{-\infty\}$
such that $f|_E\equiv -\infty$}
subset. Let $D''$ be an open set in a
   complex  manifold $X$. If $f$ is a meromorphic function on $D'\times D''$
such that $f(g,\cdot)$ extends to a meromorphic function on $X$
   for  all $g\in E$, then $f$ extends to a meromorphic function in
   a neighborhood of $E\times X\subset D'\times X$.
\end{Prop}

\begin{Lemma}
	The map $\tilde\phi_i\colon \hat G\times\C^n\to \C$ is real analytic.
\end{Lemma}

	{\bf Proof.} The question of real analyticity of
   $\tilde\phi_i$ is  local with respect to $\hat G$ so we can take a
   real  analytic coordinate neighborhood $E$  in $\hat G$, regarded as
   an  open subset of $\R$ . The map $\phi_i$ is real analytic
   in  $E\times D$ and extends therefore to a holomorphic function
   in  a  neighborhood  $D'\times D''$  of  $E\times D''$  in  the complex
   manifold  $\C\times X$.  Here we must replace $D$  by a bit smaller
   neighborhood $D''\subset D$.
      The set $E$, being an open subset of $\R$ , is nonpluripolar.
   By Proposition~\ref{Kaz}, $\phi_i$ extends to a
   meromorphic functions in a neighborhood of $E\times X$.
The restriction $\tilde\phi_i$ is therefore real analytic. \hfill $\square$

		According to the
   construction  of  Chow  scheme (see Shafarevich, \cite{S},p.65)
	every graph $\Gamma_{\tilde\phi_i}\subset\C^{2n}\subset\P^{2n}$ has
   its Chow coordinate in the Chow scheme C.  The  Chow
   coordinates  yield  a continuous mapping $f\colon \hat G \to C$. This is
   the  universal property of the Chow scheme. It follows from
   the  Theorems  of  D.  Barlet on universality of the Barlet
   space  and  on  the  equivalence  of the latter to the Chow
   scheme  in  case  of  projective  space (see Barlet, \cite{B}).
	Since $\hat G$ has finitely many components, the image in $C$ is
	has also this property. But the degree of variety is constant on
the components of $C$. This implies that the degree of the variety in
$\P^{2n}$ associated to $g\in \hat G$ is bounded. This implies that
the degrees of defining polynomials are bounded and the statement is proven.
\hfill $\square$

\end{document}